# OPTIMAL DISCRIMINATION DESIGNS[1]

By Holger Dette and Stefanie Titoff

*Ruhr-Universität Bochum*

We consider the problem of constructing optimal designs for model discrimination between competing regression models. Various new properties of optimal designs with respect to the popular $T$-optimality criterion are derived, which in many circumstances allow an explicit determination of $T$-optimal designs. It is also demonstrated, that in nested linear models the number of support points of $T$-optimal designs is usually too small to estimate all parameters in the extended model. In many cases $T$-optimal designs are usually not unique, and in this situation we give a characterization of all $T$-optimal designs. Finally, $T$-optimal designs are compared with optimal discriminating designs with respect to alternative criteria by means of a small simulation study.

**1. Introduction.** Optimal designs are frequently criticized because they are constructed from particular model assumptions before the data can be collected. Often there exist several plausible models which may be appropriate for a fit to the data. Therefore, in many applications, the data is first used to identify an appropriate model from a class of competing models and in a second step the same data is analyzed with the identified model. While the optimal design problem for the latter task has been considered by numerous authors (see, e.g., the monographs of Silvey [32], Pázman [26], Atkinson and Donev [2] or Pukelsheim [27]), much less attention has been paid to the problem of designing experiments for model discrimination. Early work was done by Stigler [34] and Studden [38], who determined optimal designs for discriminating between two nested univariate polynomials. The corresponding optimal design is called $D_s$-optimal design and minimizes the volume of

Received March 2008; revised June 2008.

[1]Supported by the Sonderforschungsbereich 475, Komplexitätsreduktion in multivariaten Datenstrukturen (Teilprojekt A2) and in part by a NIH Grant award IR01GM072876:01A1 and a BMBF-grant SKAVOE.

the confidence ellipsoid for the parameters corresponding to the extension of the smaller model. This criterion directly refers to a likelihood ratio test and was discussed by numerous authors (see, e.g., Spruill [33], Dette [10], Dette and Haller [12] or Song and Wong [35], among others). Atkinson and Fedorov [3, 4] proposed an alternative criterion, which determines a design such that the sum of squares for a lack of fit test is large. This optimality criterion is meanwhile called $T$-criterion in the statistical literature and has been considered by several authors, mostly in the context of regression models (see, e.g., Ucinski and Bogacka [39], López-Fidalgo, Tommasi and Trandafir [24] or Waterhouse et al. [40] for some recent references). The $D_s$- and $T$-optimality criteria have been studied separately without exploring the differences between both philosophies of constructing optimal designs for model discrimination.

The present paper makes an attempt to explore some relations between the—on a first glance—rather different concepts of constructing discrimination designs. In Section 2 we discuss some new properties of $T$-optimal designs and relate the $T$-optimal design problem to a problem of nonlinear approximation theory. In general, $T$-optimal designs are not unique, and in such cases we present an explicit characterization of the class of all $T$-optimal designs. In Section 3 the special case is considered where one of the competing models is linear, and here it turns out that $T$-optimal designs are in fact $D_1$-optimal (in the sense of Stigler [34]) in an extended linear regression model. This relation is then used to derive several new properties of $T$-optimal designs, especially bounds on the number of support points. In particular, it is demonstrated that in many cases the $T$-criterion yields designs which cannot be used to estimate all parameters in the extended model. Section 4 gives some more insight into the case of nonlinear regression models and also contains an extension of the results to $T$-optimality-type criteria, which are based on the Kullback-Leibler distance and have recently been proposed by López-Fidalgo, and Tommasi and Trandafir [24]. Finally, in Section 5 several examples are presented to illustrate the theoretical results. In particular, the mean squared error of parameter estimates and the power of tests based on $T$- and $D_s$-optimal designs are investigated by means of a simulation study.

**2. New properties of $T$-optimal designs.** We consider the common nonlinear regression model

$$Y = \eta(x, \theta) + \varepsilon, \tag{2.1}$$

where $\theta \in \Theta \subset \mathbb{R}^m$ is the vector of unknown parameters, and different observations are assumed to be independent. The errors are normally distributed with mean 0 and variance $\sigma^2$. In (2.1) the variable $x$ denotes the explanatory variable, which varies in the design space $\mathcal{X}$ (a more general situation with nonnormal, heteroscedastic errors is discussed in Section 4.2). We assume that $\eta$ is a continuous and real-valued function of both arguments $(x,\theta) \in \mathcal{X} \times \Theta$ and a design is defined as a probability measure $\xi$ on $\mathcal{X}$ with finite support (see Kiefer [21]). If the design $\xi$ has masses $w_i$ at the point $x_i$ ($i = 1, \ldots, k$) and $n$ observations can be made by the experimenter, this means that the quantities $w_i n$ are rounded to integers, say $n_i$, satisfying $\sum_{i=1}^{k} n_i = n$, and the experimenter takes $n_i$ observations at each location $x_i$ ($i = 1, \ldots, k$). There are numerous criteria to discriminate between competing designs, if parameter estimation in a given model is the main objective for the construction of the design (see Silvey [32], Pázman [26] or Pukelsheim [27], among others), but much less attention has been paid to the problem of developing optimal designs for model discrimination. Early work was done by Hunter and Reiner [17], Box and Hill [5] and Stigler [34]. A review on discrimination designs can be found in Hill [18]. Stigler [34] proposed a $D_s$-criterion for discriminating between two competing (nested) models. Roughly speaking, the $D_s$-optimal design yields small variances of the parameter estimates in an "extended" model. To be precise, consider the case of two rival models for the mean effect in the nonlinear regression model (2.1), say $\eta_1(x, \theta_{(1)})$ and $\eta_2(x, \theta_{(2)})$ with $\theta_{(j)} \in \Theta_{(j)} \subset \mathbb{R}^{m_j}$ ($m_j \in \mathbb{N}$, $j = 1, 2$). We assume the model $\eta_1(x, \theta_{(1)})$ is an extension of the model $\eta_2(x, \theta_{(2)})$. In other words, if the last $m_0 = m_1 - m_2$ components of the vector $\theta_{(1)} = (\theta_{(2)}, \theta_{(0)})$ vanish we obtain the model $\eta_2$, that is, $\eta_1(x, (\theta_{(2)}^T, 0^T)^T) = \eta_2(x, \theta_{(2)})$, where $0$ denotes the $(m_1 - m_2)$-dimensional vector with all components identical $0$. The $D_{m_1-m_2}$-optimality criterion is defined by the expression

$$\Phi_{D_{m_1-m_2}}(\xi) = \frac{|M_{m_1}(\xi)|}{|M_{m_2}(\xi)|}, \tag{2.2}$$

where the matrices $M_{m_1}(\xi)$ and $M_{m_2}(\xi)$ are given by

$$M_{m_1}(\xi) = \int \frac{\partial}{\partial \theta_{(1)}} \eta_1(x, \theta_{(1)}) \left( \frac{\partial}{\partial \theta_{(1)}} \eta_1(x, \theta_{(1)}) \right)^T d\xi(x) \in \mathbb{R}^{m_1 \times m_1},$$

$$M_{m_2}(\xi) = \int \frac{\partial}{\partial \theta_{(2)}} \eta_2(x, \theta_{(2)}) \left( \frac{\partial}{\partial \theta_{(2)}} \eta_2(x, \theta_{(2)}) \right)^T d\xi(x) \in \mathbb{R}^{m_2 \times m_2},$$

respectively. A $D_{m_1-m_2}$-optimal design maximizes the function $\Phi_{D_{m_1-m_2}}$ in the class of all designs, satisfying Range$(K) \subset$ Range$(M_{m_1}(\xi))$, where the matrix $K$ is defined by $K^T = (0, I_{m_1-m_2}) \in \mathbb{R}^{(m_1-m_2) \times m_1}$, $I_{m_1-m_2} \in \mathbb{R}^{(m_1-m_2) \times (m_1-m_2)}$ is the identity matrix and $0$ denotes the $(m_1 - m_2) \times m_2$ matrix with all entries identical $0$. The criterion is motivated by the likelihood ratio test for the hypothesis

$$H_0 : K^T \theta_{(1)} = 0. \tag{2.3}$$

Because the volume of the confidence ellipsoid for the parameter $K^T\theta_{(1)}$ is minimized if the function $\Phi_{D_{m_1-m_2}}(\xi)$ is maximized with respect to $\xi$ (see Pukelsheim [27]), we expect that a $D_{m_1-m_2}$-optimal design yields good power for the test of the hypothesis (2.3). The $T$-optimality criterion was introduced by Atkinson and Fedorov [3, 4], as a criterion which directly reflects the goal of model discrimination in the design of experiment and has found considerable interest in the recent literature (see, e.g., Ucinski and Bogacka [39], López-Fidalgo, Tommasi and Trandafir [24] or Waterhouse et al. [40], among many others). It does not necessarily refer to nested models and assumes that one model, say $\eta = \eta_1$ is fixed. The $T$-optimality criterion determines the design $\xi$ such that the expression

$$(2.4) \qquad \Delta(\xi) = \inf_{\theta_{(2)} \in \Theta_{(2)}} \int_{\mathcal{X}} (\eta(x) - \eta_2(x, \theta_{(2)}))^2 \, d\xi(x)$$

is maximal. The statistical interpretation of the $T$-optimality criterion is as follows. Assume that we are interested in the problem of testing the hypothesis $H_0: \eta = \eta_1$ versus $H_1: \eta = \eta_2$, which corresponds in the context of nested models to the hypotheses

$$(2.5) \qquad H_0: \theta_{(1)} = \begin{pmatrix} \theta_{(2)} \\ 0 \end{pmatrix} \quad \text{versus} \quad H_1: \theta_{(1)} \neq \begin{pmatrix} \theta_{(2)} \\ 0 \end{pmatrix}.$$

Under local alternatives of the form $\theta_{(1),n} = \begin{pmatrix} \theta_{(2)} \\ \theta_{(0)}/\sqrt{n} \end{pmatrix}$ it follows that the noncentrality parameter of the corresponding likelihood ratio test up to the factor $\sigma^2$ is given by

$$\delta^2 = \theta_{(0)}^T M_{11.2}(\xi) \theta_{(0)},$$

where $M_{11.2}(\xi)$ denotes the Schur complement of the matrix $M_{m_2}(\xi)$ in $M_{m_1}(\xi)$ and a straightforward calculation shows that

$$\delta^2 = \Delta(\xi) + o(1),$$

where the function $\eta$ in (2.4) is given by $\eta(\cdot) = \eta_1(\cdot, (\theta_{(2)}^T, \theta_{(0)}^T)^T)$. Thus a $T$-optimal design maximizes the power of the likelihood ratio test with respect to local alternatives.

The $L^2$-distance in (2.4) corresponds to the assumption of a normal distributed, homoscedastic error and alternative metrics could be used reflecting different assumptions regarding the error distribution and variance structure. For example, recently López-Fidalgo, Tommasi and Trandafir [24] proposed a Kullback-Leibler distance, which corresponds to the likelihood ratio test for the hypothesis $H_0: \eta_1 = \eta_2$ versus $H_1: \eta_1 \neq \eta_2$ under different distributional assumptions. In the present paper we will restrict ourselves to the criteria (2.2) and (2.4), but mention possible extensions of our results in the second part of Section 4.

Note that the $T$-optimality criterion, and in the case of nonlinear regression models also the $D_s$-optimality criterion, depends on the unknown parameter $\theta_{(1)}$, which may be difficult to choose in concrete applications. However, a robust version of the two optimality criteria can easily be obtained applying a sequential, Bayesian or (standardized) maximin approach (see, e.g., Atkinson and Fedorov [3], Müller and Pázman [25], Dette and Neugebauer [13, 14] or Dette [11], among many others).

For the following discussion consider the kernel

$$(2.6) \qquad \Delta(\theta_{(2)}, \xi) = \int_{\mathcal{X}} (\eta(x) - \eta_2(x, \theta_{(2)}))^2 \, d\xi(x)$$

and define for a continuous (real-valued) function $f$ on the design space $\mathcal{X}$ its sup-norm by $\|f\|_\infty = \sup_{x \in \mathcal{X}} |f(x)|$. Throughout this paper it is assumed that the infimum in (2.6) is attained for some $\theta^*_{(2)} \in \Theta_{(2)}$ and that a $T$-optimal design exists. Moreover, we assume that the regression functions $\eta_1$ and $\eta_2$ are differentiable with respect to the second argument. Our first result characterizes a $T$-optimal design as the solution of a nonlinear approximation problem.

THEOREM 2.1.
$$\sup_\xi \Delta(\xi) = \sup_\xi \inf_{\theta_{(2)} \in \Theta_{(2)}} \Delta(\theta_{(2)}, \xi) = \inf_{\theta_{(2)} \in \Theta_{(2)}} \|\eta - \eta_2(\cdot, \theta_{(2)})\|^2_\infty.$$

Moreover, if $\xi^*$ denotes a $T$-optimal design and $\theta^*_{(2)}$ is any value corresponding to the minimum of $\Delta(\theta_{(2)}, \xi^*)$ with respect to $\theta_{(2)} \in \Theta_{(2)}$, then $\theta^*_{(2)}$ corresponds to a best uniform approximation of $\eta$ by the functions $\eta(\cdot, \theta_{(2)})$, that is,

$$\inf_{\theta_{(2)} \in \Theta_{(2)}} \|\eta - \eta_2(\cdot, \theta_{(2)})\|_\infty = \|\eta - \eta_2(\cdot, \theta^*_{(2)})\|_\infty,$$

$\Delta(\xi^*) = \|\eta - \eta_2(\cdot, \theta^*_{(2)})\|^2_\infty$ and

$$(2.7) \quad \mathrm{supp}(\xi^*) \subseteq \mathcal{A} := \{x \in \mathcal{X} \mid |\eta(x) - \eta_2(x, \theta^*_{(2)})| = \|\eta - \eta_2(\cdot, \theta^*_{(2)})\|_\infty\}.$$

PROOF. A straightforward calculation shows that

$$\sup_\xi \Delta(\xi) = \sup_\xi \inf_{\theta_{(2)} \in \Theta_{(2)}} \Delta(\theta_{(2)}, \xi)$$

$$\leq \inf_{\theta_{(2)} \in \Theta_{(2)}} \sup_{x \in \mathcal{X}} |\eta(x) - \eta_2(x, \theta_{(2)})|^2 = \inf_{\theta_{(2)} \in \Theta_{(2)}} \|\eta - \eta_2(\cdot, \theta_{(2)})\|^2_\infty.$$

On the other hand, $\theta^*_{(2)}$ minimizes the function defined by (2.6) with $\xi = \xi^*$ in the set $\Theta_{(2)}$ and therefore we obtain from the equivalence theorem for $T$-optimality (see, e.g., Atkinson and Fedorov [3])

$$\inf_{\theta_{(2)} \in \Theta_{(2)}} \|\eta - \eta_2(\cdot, \theta_{(2)})\|^2_\infty \leq \|\eta - \eta_2(\cdot, \theta^*_{(2)})\|^2_\infty = \Delta(\xi^*) = \sup_\xi \Delta(\xi),$$

which proves the first assertion of the theorem. For a proof of the second part assume that the design $\xi^*$ is a $T$-optimal design and that $\theta_{(2)}^*$ minimizes the function $\Delta(\theta_{(2)}, \xi^*)$, then the function $|\eta(x) - \eta_2(x, \theta_{(2)}^*)|$ attains its maximum at any support point of $\xi^*$ (see Atkinson and Fedorov [3]) and $\theta_{(2)}^*$ corresponds to a best uniform approximation of the function $\eta$ by functions of the form $\eta_2(\cdot, \theta_{(2)})$. Therefore, the assertion follows. □

Theorem 2.1 links the $T$-optimal design problem to a problem in nonlinear approximation theory, which will be further discussed in Sections 3 and 4. Note that the theorem provides a saddle point property of the point $(\theta_{(2)}^*, \xi^*)$ although the kernel $\Delta(\theta_{(2)}, \xi)$ is in general not convex as a function of $\theta_{(2)}$. The result is particularly useful, if the best uniform approximation of the function $\eta$ by functions of the form $\eta_2(\cdot, \theta_{(2)})$ is unique, say $\eta_2(\cdot, \overline{\theta}_{(2)})$. In this case, the set $\mathcal{A}$ in (2.7) is independent of the design $\xi^*$ and the following result allows us to characterize all $T$-optimal designs.

THEOREM 2.2. *Assume that the parameter $\overline{\theta}_{(2)}$ corresponding to the best uniform approximation of the function $\eta$ by functions of the form $\eta_2(\cdot, \theta_{(2)})$ is unique and an interior point of the set $\Theta_{(2)}$.*

(a) *If a design $\xi^*$ is $T$-optimal, then*

$$(2.8) \quad \int_{\mathcal{A}} (\eta(x) - \eta_2(x, \overline{\theta}_{(2)})) \frac{\partial}{\partial \theta_{(2)}} \eta_2(x, \theta_{(2)}) \Big|_{\theta_{(2)} = \overline{\theta}_{(2)}} d\xi^*(x) = 0.$$

(b) *Conversely, assume that a design $\xi^*$ satisfies (2.8), $\mathrm{supp}(\xi^*) \subset \mathcal{A}$ and that the minimum of the function*

$$(2.9) \quad \theta_{(2)} \longrightarrow \int_{\mathcal{A}} (\eta(x) - \eta_2(x, \theta_{(2)}))^2 \, d\xi^*(x)$$

*is attained at a unique point in the interior of $\Theta_{(2)}$, then the design $\xi^*$ is $T$-optimal.*

PROOF. For a proof of part (a) we note that by Theorem 2.1 we have $\theta_{(2)}^* = \overline{\theta}_{(2)}$, $\mathrm{supp}(\xi^*) \subset \mathcal{A}$ for any $T$-optimal design $\xi^*$. Consequently, we obtain

$$\Delta(\xi^*) = \inf_{\theta_{(2)} \in \Theta_{(2)}} \int_{\mathcal{A}} (\eta(x) - \eta_2(x, \theta_{(2)}))^2 \, d\xi^*(x)$$

and the assertion follows because $\theta_{(2)}^* = \overline{\theta}_{(2)}$ corresponds to the (unique) minimum of the function on the right-hand side.

For a proof of part (b) assume that $\mathrm{supp}(\xi^*) \subset \mathcal{A}$, then it follows from Theorem 2.1

$$\sup_\xi \Delta(\xi) = \|\eta - \eta_2(\cdot, \overline{\theta}_{(2)})\|_\infty^2 = \int_\mathcal{A} (\eta(x) - \eta_2(x, \overline{\theta}_{(2)}))^2 \, d\xi^*(x)$$

$$= \inf_{\theta_{(2)} \in \Theta_{(2)}} \int_\mathcal{X} (\eta(x) - \eta_2(x, \theta_{(2)}))^2 \, d\xi^*(x)$$

because the parameter $\overline{\theta}_{(2)}$ corresponds to the unique minimum of the function (2.9). □

Roughly speaking Theorem 2.2 provides a characterization of all $T$-optimal designs by a system of linear equations, if the parameter $\bar{\theta}_{(2)}$ corresponding to the best approximation is unique, an interior point of the set $\Theta_{(2)}$ and if the cardinality of the set $\mathcal{A}$ defined in (2.7) is finite. If $\bar{\theta}_{(2)}$ is a boundary point of $\Theta_{(2)}$ an extension of condition (2.8) can easily be derived using Lagrangian multipliers.

In many applications the best uniform approximation of the function $\eta$ by functions of the form $\eta_2(\cdot, \theta_{(2)})$ is in fact unique, and sufficient conditions for this property can be found in the books of Rice [30] or Braess [7]. Note, there is an additional assumption in part (b) of Theorem 2.2 concerning the minimum of the function defined in (2.9). The answer to the question if this assumption is satisfied depends on the function $\eta$ and the parameter set $\Theta_{(2)} \subset \mathbb{R}^{m_2}$. For example, in the linear case, that is $\eta_2(x, \theta_{(2)}) = \theta_{(2)}^T f(x)$ [for an appropriate vector of regression functions $f(x)$], this assumption is always satisfied, because the Hesse-matrix of $\Delta(\theta_{(2)}, \xi)$ with respect to the parameter $\theta_{(2)}$ is given by

$$\frac{\partial^2}{\partial^2 \theta_{(2)}} \Delta(\theta_{(2)}, \xi) = 2 \cdot \int_\mathcal{X} f(x) f^T(x) \, d\xi(x),$$

and therefore positive definite, if the design $\xi$ has more than $m_2$ support points.

An exchange type algorithm for the computation of $T$-optimal designs was proposed by Atkinson and Fedorov [3]. Theorem 2.2 suggests an alternative method to determine $T$-optimal designs. In a first step the best uniform approximation of the function $\eta$ by functions of the form $\eta_2(\cdot, \theta_{(2)})$ is determined. For this calculation the Remes exchange algorithm could be used in many cases, which is a common tool in approximation theory (see Rice [30], Vol. 1, pages 171–180). The algorithm also yields the set of all possible support points $\mathcal{A}$ defined in (2.7) of $T$-optimal designs and will be illustrated in the following example. Secondly, the system of equations in (2.8) is solved to characterize all $T$-optimal designs. In contrast to the method proposed by Atkinson and Fedorov [3], this approach yields all $T$-optimal designs.

EXAMPLE 2.3. Consider the $T$-optimal design problem on the interval $[-1, 1]$ for the functions

(2.10) $\quad \eta(x) = \eta_1(x, \theta_{(1)}) = 1 + x + x^3 \quad \text{and} \quad \eta_2(x, \theta_{(2)}) = \theta_{(2)1} + \theta_{(2)2} x.$

It can be shown that the best approximation of the cubic polynomial $\eta$ by linear functions $\eta_2$ alternates at most 4 times. The Remes algorithm starts with an initial guess for the best approximation of $\eta$, say $\eta_2(\cdot, \theta_{(2)}^{(0)})$. Given an approximation $\eta_2(\cdot, \theta_{(2)}^{(k)})$ in the $k$th step one determines 4 points $x_1^{(k+1)} < \cdots < x_4^{(k+1)} \in [-1, 1]$ such that

(2.11) $\quad (\eta(x_j^{(k+1)}) - \eta_2(x_j^{(k+1)}, \theta_{(2)}^{(k)}))(\eta(x_{j+1}^{(k+1)}) - \eta_2(x_{j+1}^{(k+1)}, \theta_{(2)}^{(k)})) < 0$

$j = 1, 2, 3$ [which means that the difference $\eta(x) - \eta_2(x, \theta_{(2)}^{(k)})$ has opposite sign at the adjacent points $x_j^{(k+1)}$],

(2.12) $\quad \max_{j=1}^{4} |\eta(x_j^{(k+1)}) - \eta_2(x_j^{(k+1)}, \theta_{(2)}^{(k)})| = \|\eta - \eta_2(\cdot, \theta_{(2)}^{(k)})\|_\infty$

[at one of the points $x_j^{(k+1)}$ the function $\eta - \eta_2(\cdot, \theta_{(2)}^{(k)})$ attains its sup-norm] and

(2.13) $\quad \min_{j=1}^{4} |\eta(x_j^{(k+1)}) - \eta_2(x_j^{(k+1)}, \theta_{(2)}^{(k)})| \geq \max_{j=1}^{4} |\eta(x_j^{(k)}) - \eta_2(x_j^{(k)}, \theta_{(2)}^{(k)})|.$

In the next step the parameter $\theta_{(2)}^{(k+1)}$ is determined such that

$$\max_{j=1}^{4} |\eta(x_j^{(k+1)}) - \eta_2(x_j^{(k+1)}, \theta_{(2)}^{(k+1)})|$$

is minimal [in other words, the best approximation of the function $\eta$ by $\eta_2(\cdot, \theta_{(2)})$ with respect to the sup-norm on the set $\{x_1^{(k+1)}, \ldots, x_4^{(k+1)}\}$ is determined]. It can be shown that it is always possible to choose the points $\{x_1^{(k+1)}, \ldots, x_4^{(k+1)}\}$ such that (2.13) is satisfied (see Rice [30] and note that it is easy to satisfy (2.11) and (2.12)). We have illustrated the performance of the algorithm for the models in (2.10) in Table 1 and Figure 1, where we show the parameter $\theta_{(2)}^{(k)} = (\theta_{(2)1}^{(k)}, \theta_{(2)2}^{(k)})$, the set $\{x_1^{(k)}, \ldots, x_4^{(k)}\}$ and the approximations $\eta - \eta_2(\cdot, \theta_{(2)}^{(k)})$. Note that the algorithm stops after a few iterations which is rather typical for many examples. The algorithm yields that the best approximation is given by

$$\eta(x) - \eta_2(x, \theta_{(2)}^*) = x^3 - \tfrac{3}{4} x,$$

which yields $\mathcal{A} = \{-1, -\tfrac{1}{2}, \tfrac{1}{2}, 1\}$ for the set defined in (2.7). Because all assumptions of Theorem 2.2 are satisfied (note that the regression model $\eta_2$ is linear), the system of equations (2.8) characterizes all $T$-optimal designs. A straightforward calculation shows that the set of all $T$-optimal designs is given by the one-parametric class

$$(2.14) \qquad \xi_p^* = \begin{pmatrix} -1 & -\frac{1}{2} & \frac{1}{2} & 1 \\ -\frac{1}{6} + p & p & \frac{2}{3} - p & \frac{1}{2} - p \end{pmatrix},$$

where $p \in [\frac{1}{6}, \frac{1}{2}]$. The parameter $p$ could be chosen such that a further optimality criterion (e.g., $D$-optimality for the cubic model) is maximized in the class of all $T$-optimal designs. We finally note that the exchange type algorithm proposed by Atkinson and Fedorov [3, 4] only yields the three-point design $\xi_{1/6}^*$ as $T$-optimal design with a singular information matrix in the cubic regression model.

TABLE 1
*The iterations of the Remes algorithm for the calculation of the best approximation of the function $1 + x + x^3$ by linear polynomials $\theta_{(2)1} + \theta_{(2)2} x$*

| $k$ | $\theta_{(2)1}^{(k)}$ | $\theta_{(2)2}^{(k)}$ | $x_1^{(k)}$ | $x_2^{(k)}$ | $x_3^{(k)}$ | $x_4^{(k)}$ |
|---|---|---|---|---|---|---|
| 0 | 0.994 | 1.075 | $-0.9$ | $-0.2$ | 0.2 | 0.8 |
| 1 | 1.0000 | 1.8705 | $-1.000$ | $-0.153$ | 0.153 | 1.000 |
| 2 | 1.0000 | 1.7514 | $-1.000$ | $-0.538$ | 0.538 | 1.000 |
| 3 | 1.0000 | 1.7500 | $-1.000$ | $-0.500$ | 0.500 | 1.000 |

REMARK 2.4. It is worthwhile to mention that Theorems 2.1 and 2.2 do not require the assumption of nested models. This assumption is only needed for the statistical interpretation of the $T$- and $D_s$-optimality criterion.

**3. $D_1$- and $T$-optimal designs in linear regression models.** In this section we restrict ourselves to the case, where the regression model $\eta_2$ is a linear model, that is,

$$(3.1) \qquad \eta_2(x, \theta_{(2)}) = \theta_{(2)}^T f(x),$$

with $\theta_{(2)} \in \Theta_{(2)} = \mathbb{R}^{m_2}$. Note that the model $\eta = \eta_1$ is not necessarily linear (this case will be discussed later in this section). Moreover, the two models are not necessarily nested, except if it is stated explicitly in the following discussion. It turns out that in this case the $T$-optimal design is in fact also $D_1$-optimal in the sense of Stigler [34] for the regression model

$$(3.2) \qquad y = \theta_{(2)}^T f(x) + \beta \eta(x) + \varepsilon.$$

For a proof of this property let $\tilde{f}(x) = (f^T(x), \eta(x))^T \in \mathbb{R}^{m_2+1}$ denote the vector of regression functions in the linear regression model (3.2), let $e_{m_2+1} =$

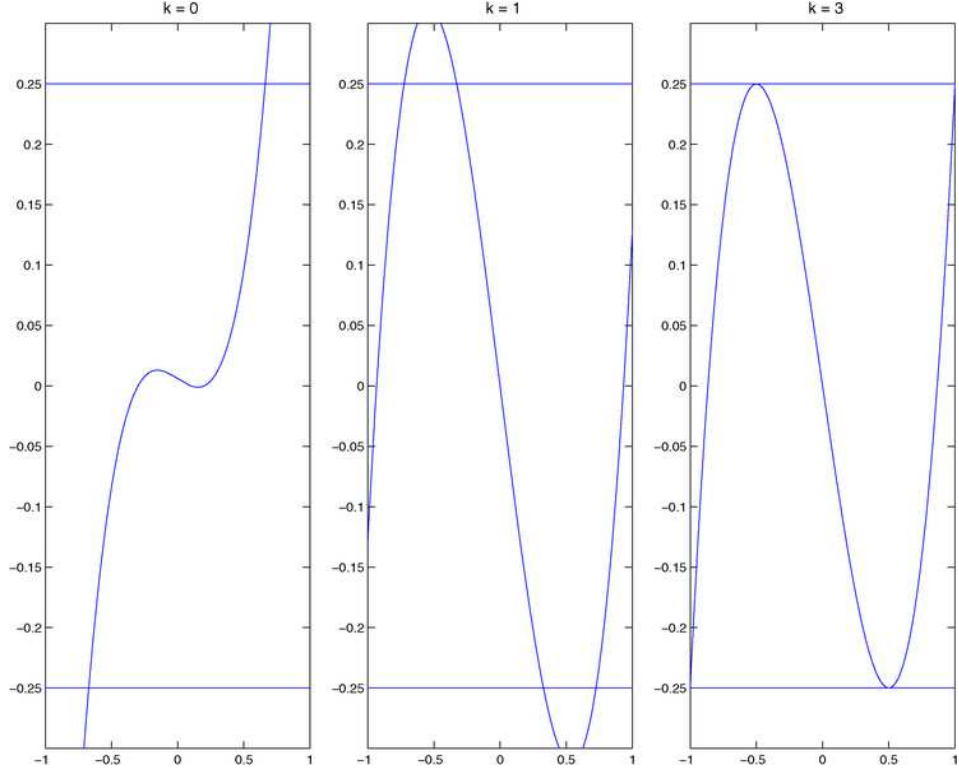

FIG. 1. *Different iteration steps of the function $1 + x + x^3 - \theta_{(2)1}^{(k)} - \theta_{(2)2}^{(k)} x$ generated by the Remes algorithm. Left panel $k = 0$, middle panel $k = 1$, right panel $k = 3$.*

$(0, \ldots, 0, 1)^T \in \mathbb{R}^{m_2+1}$ be the $(m_2 + 1)$th unit vector and define

$$(3.3) \qquad M(\xi) = \int_{\mathcal{X}} f(x) f^T(x) \, d\xi(x),$$

$$(3.4) \qquad \tilde{M}(\xi) = \int_{\mathcal{X}} \tilde{f}(x) \tilde{f}^T(x) \, d\xi(x)$$

as the information matrices in the regression model $\eta_2$ and the extended model (3.2), respectively. Recall that a $D_1$-optimal design in the regression model (3.2) satisfies $e_{m_2+1} \in \operatorname{Range}(\tilde{M}(\xi))$ and maximizes the expression

$$(e_{m_2+1}^T \tilde{M}^-(\xi) e_{m_2+1})^{-1} = \frac{\det \tilde{M}(\xi)}{\det M(\xi)}$$

(see, e.g., Stigler [34] or Studden [38]). The $D_1$-optimality criterion is a special case of the $c$-optimality criterion, which determines for a given vector $c \in \mathbb{R}^{m_2+1}$ the design $\xi$ such that the expression $(c^T \tilde{M}^-(\xi) c)^{-1}$ is maximal and the condition $c \in \operatorname{Range}(\tilde{M}(\xi))$ is satisfied (see Pukelsheim [27]). Note also that the expression $c^T \tilde{M}^-(\xi)c$ is approximately proportional to the variance of the least squares estimate of $(\theta_{(2)}^T, \beta)c$ in the regression model (3.2) (see Pukelsheim [27]). Therefore, a $D_1$-optimal design minimizes the variance of the least squares estimate of the coefficient $\beta$ in the extended regression model (3.2).

THEOREM 3.1. *Assume that (3.1) is satisfied, then a design $\xi^*$ is $T$-optimal if and only if it is $D_1$-optimal in the extended regression model (3.2).*

PROOF. Let $f(x) = (f_{(2)1}(x), \ldots, f_{(2)m_2}(x))^T$ denote the vector of functions corresponding to the first part in the linear model (3.2) and define for continuous functions $g_1, \ldots, g_k$ ($k \in \mathbb{N}$) with domain $\mathcal{X}$ the Gram determinant by

$$G(g_1, \ldots, g_k) := \left| \left( \int_{\mathcal{X}} g_i(x) g_j(x) \, d\xi(x) \right)_{i,j=1}^k \right|.$$

Then a standard result from Hilbert space theory (see Achiezer [1], page 16) shows that

$$\Delta(\xi) = \frac{G(\eta, f_{(2)1}, f_{(2)2}, \ldots, f_{(2)m_2})}{G(f_{(2)1}, f_{(2)2}, \ldots, f_{(2)m_2})} = \frac{\det \tilde{M}(\xi)}{\det M(\xi)},$$

which proves the assertion. □

In the case where the model $\eta_1(\cdot, \theta_{(1)})$ is also linear, an alternative representation for the criterion $\Delta(\xi)$ was given in Section 4.2 of Atkinson and Fedorov [3]. Theorem 3.1 provides a different interpretation of the $T$-optimality criterion and does not require the assumption of a linear model $\eta_1(\cdot, \theta_{(1)})$. In the following we derive several important conclusions from Theorem 3.1. We begin with a general result on the number of support points of $T$-optimal designs, which is a direct consequence of Corollary 8.3 in Pukelsheim [27]. Roughly speaking the number of support points of the $T$-optimal design is at most $m_2 + 1$, independently of the dimension $m_1$ of the parameter $\theta_{(1)}$ corresponding to the model $\eta_1(\cdot, \theta_{(1)})$.

COROLLARY 3.2. *Assume that (3.1) is satisfied, then there exists a $T$-optimal design $\xi^*$ with $m_2 + 1$ support points.*

We now present a refinement of this result in the case, where the design space is an interval, say $I \subset \mathbb{R}$ and the regression functions in model (3.2) form a Chebyshev system (see Karlin and Studden [20]). In many cases (with a minor additional assumption) the $T$-optimal design is supported at precisely $m_2 + 1$ well defined points, which correspond to the system under consideration and can be found explicitly. To be precise recall that a set of $k$ functions $h_1, \ldots, h_k : I \to \mathbb{R}$ is called a weak Chebyshev system (on the interval $I$) if there exists an $\varepsilon \in \{-1, 1\}$ such that the inequality

$$
(3.5) \qquad \varepsilon \cdot \begin{vmatrix} h_1(x_1) & \ldots & h_1(x_k) \\ \vdots & \ddots & \vdots \\ h_k(x_1) & \ldots & h_k(x_k) \end{vmatrix} \geq 0
$$

holds for all $x_1, \ldots, x_k \in I$ with $x_1 < x_2 < \cdots < x_k$. If the inequality in (3.5) is strict, then $\{h_1, \ldots, h_k\}$ is called a Chebyshev system. It is well known (see Karlin and Studden [20], Theorem II 10.2) that if $\{h_1, \ldots, h_k\}$ is a Chebyshev system, then there exists a unique function, say $\sum_{i=1}^k c_i^* h_i(x) = c^{*T} h(x)$, $(h = (h_1, \ldots, h_k)^T)$ with the following properties

$$
\begin{aligned}
&\text{(i)} && |c^{*T} h(x)| \leq 1 \qquad \forall x \in I \\
(3.6) \quad &\text{(ii)} && \text{there exist } k \text{ points } x_1^* < \cdots < x_k^* \text{ such that} \\
& && c^{*T} h(x_i^*) = (-1)^i, \qquad i = 1, \ldots, k.
\end{aligned}
$$

The function $c^{*T} h(x)$ is called Chebyshev polynomial, and we say that it is alternating at the points $x_1^*, \ldots, x_k^*$. The points $x_1^*, \ldots, x_k^*$ are called Chebyshev points and need not to be unique. They are unique in most applications, in particular if $1 \in \text{span}\{h_1, \ldots, h_k\}, k \geq 1$ and $I$ is a bounded and closed interval, where in this case $x_1^* = \min_{x \in I} x$, $x_k^* = \max_{x \in I} x$. It is well known (see Studden [36], Pukelsheim and Studden [28] or Imhof and Studden [19], among others) that in many cases $c$-optimal designs in regression models are supported at the Chebyshev points. The following result shows that a similar statement can be made for $T$-optimal designs.

THEOREM 3.3. *Assume that* (3.1) *is satisfied, that the design space is an interval, say $\mathcal{X} = I \subset \mathbb{R}$ and that $\{f_1, \ldots, f_{m_2}\}$ is a Chebyshev system on the interval $I$. In this case the set $\mathcal{A}$ defined in* (2.7) *has at least $m_2 + 1$ points.*

*Moreover, assume that additionally $\{f_1, \ldots, f_{m_2}, \eta\}$ is also a Chebyshev system on the interval $I$ and*

$$
\begin{vmatrix} f_1(x_1) & \ldots & f_1(x_{m_2}) & 0 \\ \vdots & \ddots & \vdots & \vdots \\ f_{m_2}(x_1) & \ldots & f_{m_2}(x_{m_2}) & 0 \\ \eta(x_1) & \ldots & \eta(x_{m_2}) & 1 \end{vmatrix} \neq 0
$$

*for all $x_1, \ldots, x_{m_2} \in I$ satisfying $x_1 < \cdots < x_{m_2}$. Let $x_1^* < \cdots < x_{m_2+1}^*$ denote $m_2 + 1$ Chebyshev points satisfying* (3.6) *and define $\xi^*$ as the design which*

*has weights*

$$w_i^* = \frac{|u_i|}{\sum_{j=1}^{m_2+1} |u_j|}$$

*at the points* $x_i^*$ $(i = 1, \ldots, m_2 + 1)$, *where* $u = (u_1, \ldots, u_{m_2+1})^T = (X^T X)^{-1} \times X^T e_{m_2+1}$, *and the matrix* $X$ *is defined by*

$$X = \begin{pmatrix} f_1(x_1^*) & \cdots & f_1(x_{m_2+1}^*) \\ \vdots & \ddots & \vdots \\ f_{m_2+1}(x_1^*) & \cdots & f_{m_2+1}(x_{m_2+1}^*) \end{pmatrix}$$

*(here we put* $f_{m_2+1} = \eta$*). Then* $\xi^*$ *is a $T$-optimal design.*

PROOF. It follows from Theorem 1.1 in Chapter IX of Karlin and Studden [20] that the best uniform approximation of the function $\eta$ by functions of the form $\eta_2(x, \theta_{(2)}) = \theta_{(2)}^T f(x)$ is unique. By Theorem 2.1 the support of a $T$-optimal design is contained in the set

$$\mathcal{A} = \left\{ x \in I \,\bigg|\, \left| \eta(x) - \sum_{j=1}^{m_2} \overline{\theta}_{(2)j} f_j(x) \right| = \left\| \eta - \sum_{j=1}^{m_2} \overline{\theta}_{(2)j} f_j \right\|_\infty \right\},$$

where the parameters $\overline{\theta}_{(2)1}, \ldots, \overline{\theta}_{(2)m_2}$ correspond to the best uniform approximation of $\eta$ by linear combinations of $f_1, \ldots, f_{m_2}$. Theorem 1.1 in Karlin and Studden [20] also shows that the cardinality of the set $\mathcal{A}$ is at least $m_2 + 1$ and the first assertion follows.

For a proof of the second part we note that by Theorem 3.1 the $T$-optimal design problem is equivalent to the $D_1$-optimal design problem in the extended regression model (3.2). Because this is exactly the $e_{m_2+1}$-optimal design problem it follows from Kiefer and Wolfowitz [22] (see also Studden [36]) that the $T$-optimal design is supported at $m_2 + 1$ points satisfying (3.6). The formula for the corresponding weights is now a direct consequence of Corollary 8.9 in Pukelsheim [27]. □

If, under the assumptions of Theorem 3.3 there exist exactly $m_2 + 1$ uniquely determined Chebyshev points, then any $T$-optimal design is supported at precisely $m_2 + 1$ points. This situation is rather typical in applications. Note that in Example 2.3 ($m_2 = 2$) the functions $\{1, x\}$ form a Chebyshev system. Thus the first part of Theorem 3.3 implies that the set $\mathcal{A}$ in (2.7) has at least cardinality 3 (in fact its cardinality is 4). On the other hand, the system $\{1, x, x^3\}$ is not a Chebyshev system on the interval $[-1, 1]$, because the polynomial $x^3 - \frac{3}{4}x$ has 3 roots in the interval $[-1, 1]$. As a consequence the second part of Theorem 3.3 is not applicable here. In fact, there exist an infinite number of $T$-optimal designs with 4 support points indicating that the Chebyshev property of the system $\{f_1, \ldots, f_{m_2}, \eta\}$ is really necessary in this context.

In the following we specialize the result of Theorem 3.1 to the case, where the model $\eta_1$ is in fact an extension of the linear regression model (3.1), that is $\theta_{(1)} = (\theta_{(2)}^T, \theta_{(0)}^T)^T$,

$$\eta(x) = \eta_1(x, \theta_{(1)}) = \theta_{(2)}^T f(x) + \theta_{(0)}^T g(x), \tag{3.7}$$

where $g(x) = (g_1(x), \ldots, g_{m_0}(x))^T$ is a further vector of regression functions and $m_0 + m_2 = m_1$. In this case, Theorem 3.1 can be slightly simplified.

COROLLARY 3.4. *Assume that* (3.1) *and* (3.7) *are satisfied, then a design $\xi^*$ is $T$-optimal if and only if it is $D_1$-optimal in the extended regression model*

$$y = \theta_{(2)}^T f(x) + \beta \phi(x) + \varepsilon, \tag{3.8}$$

*where* $\phi(x) = \theta_{(0)}^T g(x)$.

PROOF. From Theorem 3.1 and its proof it follows that a design is $T$-optimal if and only if it maximizes

$$\frac{\det \tilde{M}(\xi)}{\det M(\xi)} = \frac{G(\eta, f_{(2)1}, f_{(2)2}, \ldots, f_{(2)m_2})}{G(f_{(2)1}, f_{(2)2}, \ldots, f_{(2)m_2})} \tag{3.9}$$

$$= \frac{G(\theta_{(0)}^T g, f_{(2)1}, f_{(2)2}, \ldots, f_{(2)m_2})}{G(f_{(2)1}, f_{(2)2}, \ldots, f_{(2)m_2})},$$

where the matrix $\tilde{M}(\xi)$ is defined by (3.4). The last equality follows from (3.7) and the multi-linearity of the Gram determinant. Therefore the $T$-optimal design is $D_1$-optimal in the regression (3.8). □

We conclude this section with an alternative interpretation of the $T$-optimality criterion as a compound criterion in the situation considered in Corollary 3.4. To be precise, we define the $m_0 = m_1 - m_2$ regression models

$$y = \theta_{(2)}^T f(x) + \beta_j g_j(x) + \varepsilon, \qquad j = 1, \ldots, m_0.$$

Then, by Theorem 3.1, the $T$-optimal design for discriminating between $\eta_2$ and the $j$th model $(\theta_{(2)}^T, \beta_j) \tilde{f}_j(x)$ with $\tilde{f}_j(x) = (f^T(x), g_j(x))^T$ maximizes

$$\Delta_j(\xi) = \frac{\det \tilde{M}_j(\xi)}{\det M(\xi)} = \frac{G(g_j, f_{(2)1}, \ldots, f_{(2)m_2})}{G(f_{(2)1}, \ldots, f_{(2)m_2})}, \qquad j = 1, \ldots, m_0,$$

where

$$\tilde{M}_j(\xi) = \int_{\mathcal{X}} \tilde{f}_j(x) \tilde{f}_j^T(x) \, d\xi(x)$$

and the matrix $M(\xi)$ is defined in (3.3). The proof of the next result is now a direct consequence of the representation (3.9) and the multilinearity of the Gram determinant.

COROLLARY 3.5. *A $T$-optimal design for discriminating between the models (3.1) and (3.7) maximizes the weighted average*

$$\Delta(\xi) = \sum_{j=1}^{m_0} \theta_{(0)j} \Delta_j(\xi),$$

*where $\theta_{(0)j}$ denotes the $j$th component of the vector $\theta_{(0)}$ in (3.7).*

Note that by Corollary 3.5 the $T$-optimal design for discriminating between the models (3.1) and (3.7) can be interpretated as a compound optimality criterion in the sense of Läuter [23] and therefore results for calculating optimal designs with respect to compound criteria can be used to find $T$-optimal designs (see, e.g., Pukelsheim [27], Cook and Wong [9] or Clyde and Chaloner [8], among many others).

## 4. Further discussion.

4.1. *Some comments on nonlinear models.* As mentioned before, in general Theorems 2.1 and 2.2 link the $T$-optimal design problem to a problem in nonlinear approximation theory, which has a long history in mathematics (see Braess [7] or Rice [30]), and is substantially more difficult to analyze compared to the linear case considered in Section 3. We will now indicate how this theory can be used to transfer some of the results of Section 3 to the nonlinear case. For this we assume that the design space $\mathcal{X}$ is an interval and that the function $\eta_2$ is continuous on $\mathcal{X} \times \Theta_{(2)}$. The following definition is taken from Rice [30].

DEFINITION 4.1. *The class of functions $\mathcal{M} = \{\eta_2(\cdot, \theta_{(2)}) | \ \theta_{(2)} \in \Theta_{(2)}\}$ has property Z of degree $m = m(\theta_{(2)}^*)$ at the point $\theta_{(2)}^* \in \Theta_{(2)}$, if for any $\theta_{(2)} \in \Theta_{(2)}$ with $\theta_{(2)} \neq \theta_{(2)}^*$ the difference $\eta_2(x, \theta_{(2)}^*) - \eta_2(x, \theta_{(2)})$ has at most $m - 1$ zeros in $\mathcal{X}$.*

*The class of functions $\{\eta_2(\cdot, \theta_{(2)}) | \theta_{(2)} \in \Theta_{(2)}\}$ is called locally solvent of degree $m = m(\theta_{(2)}^*)$ at the point $\theta_{(2)}^* \in \Theta_{(2)}$, if given a set $\{x_1, \ldots, x_m\} \subset \mathcal{X}$ and $\varepsilon > 0$, there exists a number $\delta = \delta(\theta_{(2)}^*, \varepsilon, x_1, \ldots, x_m) > 0$ such that the inequalities*

$$|Y_i - \eta_2(x_i, \theta_{(2)}^*)| < \delta \qquad (i = 1, \ldots, m)$$

*imply the existence of a solution $\theta_{(2)} \in \Theta_{(2)}$ of the system of nonlinear equations*

$$\eta_2(x_i, \theta_{(2)}) = Y_i, \qquad i = 1, 2, \ldots, m$$

*which satisfies*

$$\|\eta_2(\cdot,\theta_{(2)}) - \eta_2(\cdot,\theta_{(2)}^*)\|_\infty < \varepsilon.$$

*The class $\mathcal{M}$ is called varisolvent if at each point the local solvency property and property Z are satisfied with the same degree.*

Examples of varisolvent families include sums of exponentials and rational functions (see Rice [30]). If the class of functions $\{\eta_2(\cdot,\theta_{(2)})|\ \theta_{(2)} \in \Theta_{(2)}\}$ is varisolvent, the following result gives a rough estimate of the number of support points of the $T$-optimal design. The proof can be found in Braess [7].

THEOREM 4.2. *Assume that the class of functions $\mathcal{M} = \{\eta_2(\cdot,\theta_{(2)})|\theta_{(2)} \in \Theta_{(2)}\}$ is varisolvent and that $\eta$ is a continuous function on $\mathcal{X}$ such that $\eta - \eta_2(\cdot,\bar{\theta}_{(2)})$ is not constant. The function $\eta_2(\cdot,\bar{\theta}_{(2)})$ is a best approximation of the function $\eta$ if and only if the difference $\eta - \eta_2(\cdot,\bar{\theta}_{(2)})$ alternates $m(\bar{\theta}_{(2)})+1$ times, that is, there exists at least $m(\bar{\theta}_{(2)})+1$ points $x_0^* < \cdots < x_{m(\bar{\theta}_{(2)})}^*$ in $\mathcal{X}$ such that*

$$\eta(x_i^*) - \eta_2(x_i^*,\bar{\theta}_{(2)}) = \varepsilon(-1)^i \|\eta - \eta_2(\cdot,\bar{\theta}_{(2)})\|_\infty, \qquad i = 0,\ldots,m(\bar{\theta}_{(2)}),$$

*where $\varepsilon \in \{-1,1\}$.*

Theorem 4.2 gives some hint of the number of support points of the $T$-optimal design. By this result, there exists a best approximation of $\eta$ by functions of the form $\eta_2(\cdot,\theta_{(2)})$ $(\theta_{(2)} \in \Theta_{(2)})$, such that $r^* = \eta - \eta_2(\cdot,\bar{\theta}_{(2)})$ alternates at least $m(\bar{\theta}_{(2)})+1$ times. In many cases there are no other points in $\mathcal{X}$ where the difference $r^*$ attains its maximum, and it follows from Theorem 2.1 that the $T$-optimal design has at most $m(\bar{\theta}_{(2)})+1$ support points. We illustrate this heuristic argument by an example, where we consider sums of exponentials.

EXAMPLE 4.3. Assume that

$$\eta_2(x,\theta_{(2)}) = \sum_{j=1}^{m_2} \theta_{(2)2j-1}\ e^{-\theta_{(2)2j}x},$$

where $x \in \mathcal{X} \subset [0,\infty)$, $\theta_{(2)2k-1} \in \mathbb{R}$, $\theta_{(2)2k} \in \mathbb{R}^+$ $(k=1,\ldots,m_2)$ and the design space is a compact interval. Models of this type have numerous applications in pharmacokinetics (see, e.g., Shargel and Yu [31] or Rowland [29]). It follows from Braess [7], pages 190–191, that for each

$$u(x) = \sum_{j=1}^{l} a_j e^{-b_j x}$$

with $b_1, \ldots, b_l \neq 0$, the class of functions $\mathcal{F} = \{\eta_2(\cdot, \theta_{(2)}) \mid \theta_{(2)} \in \mathbb{R}^{2m_2}, \theta_{(2)2j} \in \mathbb{R}^+; j = 1, \ldots, m_2\}$ is locally solvent at $u$ of order $m_2 + l$. Similarly, the class $\mathcal{F}$ has property $Z$ of degree $m_2 + l$ at $u$, and therefore it is varisolvent at $u$ of degree $m_2 + l$. If $\eta = \eta_1$ is a continuous function and $\eta_2(\cdot, \bar{\theta}_{(2)})$ is the best approximation of $\eta$, it follows from Theorem 4.2 that the difference $\eta - \eta_2(\cdot, \bar{\theta}_{(2)})$ alternates (at least) $m(\bar{\theta}_{(2)}) + 1 = m_2 + l(\bar{\theta}_{(2)}) + 1$, where $l(\bar{\theta}_{(2)})$ denotes the number of non-vanishing coefficients among $\bar{\theta}_{(2)1}, \bar{\theta}_{(2)3}, \ldots, \bar{\theta}_{(2)2m_2-1}$ in $\eta_2(x, \bar{\theta}_{(2)})$. By Theorem 2.1 the support points of a $T$-optimal design must be among the points, where the function $\eta - \eta_2(\cdot, \bar{\theta}_{(2)})$ attains its maximum. If none of the coefficients $\bar{\theta}_{2(2j-1)}$ vanishes, the cardinality of the set $\mathcal{A}$ in (2.7) is at least $2m_2 + 1$.

The upper bound on the cardinality of the set $\mathcal{A}$ depends on the particular properties of the function $\eta = \eta_1$ and is in many cases close to the lower bound $2m_2 + 1$. For example, if $\eta_1$ is also a sum of exponentials, say

$$\eta_1(x, \theta_{(1)}) = \sum_{j=1}^{m_1} \theta_{(1)2j-1} \ e^{-\theta_{(1)2j}x},$$

$\theta_{(1)2j-1} \in \mathbb{R}$, $\theta_{(1)2j} \in \mathbb{R}^+$, where $m_1 = m_2 + m_0 > m_2$, the difference $r^* = \eta_1 - \eta_2(\cdot, \bar{\theta}_{(2)})$ consists of at most $m_1 + m_2$ different exponential terms. Because of the Chebyshev property of the function $\{e^{a_j x} \mid j = 1, \ldots, l\}$ on the nonnegative line $(0, \infty)$ (see Karlin and Studden [20]) it follows that the derivative of the difference $r^*$ (which is also a sum of at most $m_1 + m_2$ exponential terms) has at most $m_1 + m_2 - 1$ roots. Observing that $\lim_{x \to \infty} r^*(x) = 0$ it therefore follows that there exist at most $m_1 + m_2$ alternating points of the difference $r^*$. Moreover, if the cardinality of the set $\mathcal{A}$ is exactly $m_1 + m_2$, then a boundary point of the design space $\mathcal{X}$ is an element of the set $\mathcal{A}$. Consequently any $T$-optimal design has at most $m_1 + m_2$ support points. Note that the number of parameters in the exponential models $\eta_1$ and $\eta_2$ is $2m_1$ and $2m_2$, respectively. Because $m_2 < m_1$ the $T$-optimal design cannot be used to estimate all parameters in the extended model $\eta_1$. For example, if $m_1 = m_2 + 1$, it follows from these arguments that a $T$-optimal design has precisely $2m_2 + 1$ support points, although the model $\eta_1$ has $2m_2 + 2$ parameters.

4.2. *T-optimality based on the Kullback–Leibler distance.* Recently López-Fidalgo, Tommasi and Trandafir [24] considered a generalization of the $T$-optimality criterion, which is based on the popular Kullback–Leibler (KL)-distance. The general criterion addresses the problem of a nonnormal error distribution and heteroscedasticity in model (2.1). It reduces to the $T$-criterion in the case of normal and homoscedastic data. We briefly indicate that the results of the previous sections can be easily extended to this more general class of optimality criteria.

Following López-Fidalgo, Tommasi and Trandafir [24] we specify the two different models by their densities, say $f_j(y, x, \theta_{(j)}, \sigma^2); \theta_{(j)} \in \Theta_{(j)}; j = 1, 2$, where $\sigma^2$ is a nuisance parameter corresponding to the variances of the responses. We fix one model, say $f(y, x, \sigma^2) = f_1(y, x, \theta_{(1)}, \sigma^2)$, and consider for a design $\xi$ the optimality criterion

$$\Delta_{\mathrm{KL}}(\xi) = \inf_{\theta_{(2)} \in \Theta_{(2)}} \int_{\mathcal{X}} d_{\mathrm{KL}}(f, f_2, x, \theta_{(2)}) \, d\xi(x), \tag{4.1}$$

where (for any $x \in \mathcal{X}$)

$$d_{\mathrm{KL}}(f, f_2, x, \theta_{(2)}) = \int f(y, x, \sigma^2) \log\left\{\frac{f(y, x, \sigma^2)}{f_2(y, x, \theta_{(2)}, \sigma^2)}\right\} dy$$

denotes the KL-distance between the "true" model $f$ and the alternative model $f_2(y, x, \theta_{(2)}, \sigma^2)$. A KL-optimal design $\xi^*_{\mathrm{KL}}$ maximizes $\Delta_{\mathrm{KL}}(\xi)$ in the class of all designs. The goal of this criterion is to determine designs maximizing the power of the likelihood ratio test for the hypotheses

$$H_0 : f(x, y, \sigma^2) = f_2(x, y, \theta_{(2)}, \sigma^2) \quad \text{vs.} \quad H_1 : f(y, x, \sigma^2) = f_1(y, x, \theta_{(1)}, \sigma^2)$$

for the "worst" choice $\theta_{(2)} \in \Theta_{(2)}$. Similar arguments as given in the proof of Theorem 2.1 show that

$$\sup_{\xi} \Delta_{\mathrm{KL}}(\xi) = \inf_{\theta_{(2)} \in \Theta_{(2)}} \|d_{\mathrm{KL}}(f, f_2, \cdot, \theta_{(2)})\|_{\infty} = \|d_{\mathrm{KL}}(f, f_2, \cdot, \theta^*_{(2)})\|_{\infty},$$

where $\theta^*_{(2)}$ corresponds to the minimum in (4.1) for the design $\xi^*_{\mathrm{KL}}$ and the support of a KL-optimal design $\xi^*_{\mathrm{KL}}$ satisfies

$$\mathrm{supp}(\xi^*_{\mathrm{KL}}) \subset \mathcal{A}_{\mathrm{KL}} = \{x \in \mathcal{X} | d_{\mathrm{KL}}(f, f_2, x, \theta^*_{(2)}) = \|d_{\mathrm{KL}}(f, f_2, \cdot, \theta^*_{(2)})\|_{\infty}\}.$$

This means that the KL-optimal design problem is closely related to the problem of determining the best uniform approximation of the function $\eta \equiv 0$ by the (nonlinear) parametric family

$$\{d_{\mathrm{KL}}(f, f_2, \cdot, \theta_{(2)}) \mid \theta_{(2)} \in \Theta_{(2)}\}. \tag{4.2}$$

Therefore, all results of the previous sections remain valid, where the class $\{\eta_2(\cdot, \theta_{(2)}) \mid \theta_{(2)} \in \Theta_{(2)}\}$ has to be replaced by the set defined in (4.2) and the function $\eta = \eta_1$ is given by $\eta(x) \equiv 0$. We will illustrate these ideas with an example for heteroscedastic regression models with normal distributed responses.

EXAMPLE 4.4. We consider the problem of discriminating between two regression models with heteroscedastic but normally distributed errors, that is,

$$P_j^{Y|x} \sim \mathcal{N}(\eta_j(x, \theta_{(j)}), (1-x^2)^{-1}), \qquad j = 1, 2,$$

where $\eta_1(x, \theta_{(1)}) = \eta(x) = 8x^3$ is a cubic, $\eta_2(x, \theta_{(2)}) = \theta_{(2)1} + \theta_{(2)2}x$ a linear polynomial and the explanatory variable satisfies $x \in (-1, 1)$. $D$-optimal designs for polynomial regression models with variance function $(1 - x^2)^{-1}$ have been studied extensively in the literature (see, e.g., Fedorov [16]), but discrimination designs have not been considered so far. If $f_j(y, x, \theta_{(j)})$ denotes the density of $P_j^{Y|x}$ with respect to the Lebesgue measure it follows by a straightforward but tedious calculation that

$$(4.3) \quad d_{\mathrm{KL}}(f, f_2, x, \theta_{(2)}) = (1 - x^2)(8x^3 - \theta_{(2)2}x - \theta_{(2)1})^2,$$

and the best uniform approximation of the function $\eta \equiv 0$ by functions of the form (4.3) is unique and given by $d_{\mathrm{KL}}(f, f_2, x, \bar{\theta}_{(2)}) = (8x^3 - 4x)^2(1 - x^2)$ with corresponding set

$$\mathcal{A}_{\mathrm{KL}} = \{-\tfrac{1}{2}\sqrt{2 + \sqrt{2}}, -\tfrac{1}{2}\sqrt{2 - \sqrt{2}}, \tfrac{1}{2}\sqrt{2 - \sqrt{2}}, \tfrac{1}{2}\sqrt{2 + \sqrt{2}}\}$$

and $\|d_{\mathrm{KL}}(f, f_2, x, \bar{\theta}_{(2)})\|_\infty = 1$. The analogue of Theorem 2.2 shows that all KL-optimal designs are supported in $\mathcal{A}_{\mathrm{KL}}$ and characterized by the analogue of (2.8), which yields

$$\int_{\mathcal{A}_{\mathrm{KL}}} \frac{\partial}{\partial \theta_{(2)}} d_{\mathrm{KL}}(f, f_2, x, \theta_{(2)}^*) \, d\xi^*(x) = -2 \int_{\mathcal{A}_{\mathrm{KL}}} (1 - x^2)(8x^3 - 4x) \binom{1}{x} d\xi^*(x)$$
$$= 0.$$

A straightforward calculation shows that all KL-optimal designs are given by the one-parametric class

$$\xi_{\mathrm{KL}}^* =$$
$$\begin{pmatrix} \dfrac{-\sqrt{2+\sqrt{2}}}{2} & \dfrac{-\sqrt{2-\sqrt{2}}}{2} & \dfrac{\sqrt{2-\sqrt{2}}}{2} & \dfrac{\sqrt{2+\sqrt{2}}}{2} \\ p & \dfrac{(2-\sqrt{2}) + 4p(\sqrt{2}-1)}{4} & \dfrac{\sqrt{2} - 4p(\sqrt{2}-1)}{4} & \dfrac{1}{2} - p \end{pmatrix},$$

where $p \in [0, \tfrac{1}{2}]$. We finally note that the algorithm proposed by López-Fidalgo, Tommasi and Trandafir [24] yields to the 3-point design obtained for $p = 1/2$, which cannot be used for estimating the parameters in the cubic model.

**5. Examples.** In this section we compare $T$- and $D_s$-optimal designs with respect to their power properties and estimation error by means of a simulation study. We begin with the case of discriminating between two polynomials of degree $m_2 - 1$ and $m_1 - 1 = m_0 + m_2 - 1$ on a nonnegative interval. Our second example considers a nonlinear case, namely exponential regression models.

5.1. *Polynomial regression.* Consider the polynomial regression models

$$\eta_2(x, \theta_{(2)}) = \theta_{(2)1} + \theta_{(2)2}x + \cdots + \theta_{(2)m_2}x^{m_2-1},$$

$$\eta_1(x, \theta_{(1)}) = \theta_{(1)1} + \theta_{(1)2}x + \cdots + \theta_{(1)m_2}x^{m_2-1} + \cdots + \theta_{(1)m_0+m_2}x^{m_0+m_2-1},$$

where the explanatory variable $x$ varies in a nonnegative interval, say $I \subset [0,\infty)$. Note that under the additional assumption of positive coefficients $\theta_{(1)1+m_2}, \ldots, \theta_{(1)m_0+m_2}$ the two systems of functions

$$\begin{aligned}
& \{1, x, \ldots, x^{m_2-1}, \theta_{(1)1} + \theta_{(1)2}x + \cdots + \theta_{(1)m_0+m_2}x^{m_0+m_2-1}\} \\
& \{1, x, \cdots, x^{m_2-1}\}
\end{aligned} \quad (5.1)$$

form a Chebyshev system on the interval $I$ and that the number of corresponding Chebyshev points is exactly $m_2 + 1$ (see Karlin and Studden [20], page 9). Consequently, Theorem 3.3 is applicable here and any $T$-optimal design is supported at $m_2 + 1$ points. We note that in the case $m_0 > 1$ the $T$-optimal design cannot be used for the $F$-test, which is commonly applied to discriminate between the two nested polynomials and requires at least $m_0 + m_2$ different design points. Note also that this problem was already observed by Atkinson and Donev [2] in the case $m_2 = 1$ and $m_0 = 2$ (see Example 20.2 in this reference). The results in the present paper show that this situation is not an exception but rather typical for discrimination designs constructed from the $T$-optimality criterion.

If the system in (5.1) is not a Chebyshev system the results of Sections 2 and 3 indicate that there exist several $T$-optimal designs. For example, consider the case of discriminating between a linear and a cubic polynomial, that is, $m_2 = 2$, $m_0 = 2$ on the interval $[-1, 1]$. For the cubic model we investigate the model

$$(5.2) \qquad \eta(x) = 1 + x + c_0 x^2 + d_0 x^3.$$

Some $T$-optimal designs for various values of the parameters $c_0$ and $d_0$ are given in Table 2.

The $T$-optimal design obtained from the algorithm of Atkinson and Fedorov [3] for the parameters $c_0 = 0$ and $d_0 = 1$ has weights $1/6$, $1/2$ and $1/3$ at the points $-1/2$, $1/2$ and $1$. This design corresponds to the choice $p = 1/6$ in Example 2.3 and will be called $T_{1/6}$-optimal design in this example. In order to compare the different designs with respect to their ability to discriminate between a cubic and a linear regression model by the common $F$-test we have modified the $T_{1/6}$-optimal design slightly and have put 2% of the observations at a fourth point, namely the left boundary of the design space. A further $T$-optimal design with four support points is obtained from formula (2.14) with $p = 1/3$ and denoted as $T_{1/3}$-optimal design. Stigler [34] proposed the $D_2$-criterion for the construction of a discriminating design between a linear and a cubic model. If $M_1(\xi)$ and $M_3(\xi)$ denote the information matrices of a design in the linear and cubic model, respectively, the corresponding $D_2$-optimal design maximizes $|M_3(\xi)|/|M_1(\xi)|$ and has weights 1/5, 3/10, 3/10 and 1/5 at the points $-1$, $-0.408$, 0.408 and 1 (see Studden [37]).

TABLE 2
*T-optimal designs for discriminating between a linear and a cubic polynomial given in (5.2) for a special choice of the parameters $c_0$ and $d_0$ (the parameter $z$ satisfies $z \in \mathbb{R} \setminus \{0\}$). In the case $(c_0, d_0) = (0, z)$ the T-optimal design is not uniquely determined*

| $c_0$ | $d_0$ | $x_1$ | $x_2$ | $x_3$ | $\omega_1$ | $\omega_2$ | $\omega_3$ |
|---|---|---|---|---|---|---|---|
| 0 | $z$ | $-0.5$ | 0.5 | 1 | 1/6 | 1/2 | 1/3 |
| $z$ | 0 | $-1$ | 0 | 1 | 1/4 | 1/2 | 1/4 |
| $z$ | $z$ | $-1$ | 0.33 | 1 | 1/6 | 1/2 | 1/3 |
| $z$ | $z$ | $-1$ | $-0.33$ | 1 | 1/3 | 1/2 | 1/6 |
| $2z$ | $z$ | $-1$ | 0.2 | 1 | 1/5 | 1/2 | 3/10 |
| $z$ | $2z$ | $-0.77$ | 0.411 | 1 | 1/6 | 1/2 | 1/3 |
| $-2z$ | $z$ | $-1$ | $-0.2$ | 1 | 3/10 | 1/2 | 1/5 |
| $z$ | $-2z$ | $-1$ | $-0.411$ | 0.77 | 1/3 | 1/2 | 1/6 |

We have conducted a small simulation study and generated normally distributed random variables with mean given by (5.2) and variance $\sigma^2 = 0.1$, where the design was either the $T_{1/3}$-optimal, the (modified) $T_{1/6}$-optimal or the $D_2$-optimal design. In Figure 2 we display the power function of the $F$-test for the hypothesis of a linear regression $H_0 : (c_0, d_0) = (0, 0)$ for various choices of the parameters $c_0$ and $d_0$. The level is 5% and the sample size is $n = 50$. We have considered three values for the parameter $c_0$ and display the power as a function of the parameter $d_0$. The solid line corresponds to the power function of the $F$-test based on the (modified) $T_{1/6}$-optimal design, while the dotted and dashed line refer to the $T_{1/3}$-optimal and $D_2$-optimal design, respectively. If $c_0 = 0$ the curves are almost identical if $d_0$ is also small, while we observe some advantages for the $T_{1/3}$- and $D_2$-optimal design for moderate and large values of $d_0$. Here the $T_{1/3}$-optimal design has the best performance (see the left panel in Figure 2). The case of a positive parameter $c_0 = 0.05$, $c_0 = 0.1$ corresponds to an alternative. For small values of $d_0$, the $D_2$-optimal and $T_{1/3}$-optimal design seem to have better discrimination properties than the $T_{1/6}$-optimal design, while the opposite behavior is observed if $d_0$ is large (see the middle and right panel in Figure 2). Next we consider the situation where $d_0$ is fixed and the parameter $c_0$ is varied. If $d_0 = 0$ the $D_2$-optimal design always yields more power than both $T$-optimal designs, where the $T_{1/3}$-optimal design shows some advantages compared to the $T_{1/6}$-optimal design (see the left panel in Figure 3). For larger values of $d_0$ the situation is similar. If $d_0 = 0.05$ all three designs yield very similar results for small values of the parameter $c_0$, while for larger values of $c_0$ the $T_{1/3}$- and $D_2$-optimal design yield more power than the $T_{1/6}$-optimal design. Finally, in the case $d_0 = 0.1$ the $T_{1/6}$-optimal design should be preferred for small values of $c_0$ if model discrimination is the main goal of the experimenter. On the other hand, if $d_0$ is large, the $D_2$-optimal design has the best performance and both $T$-optimal designs show a similar behavior (see the right panel in Figure 3). Summarizing these observations, we conclude that the superiority of one of the two discrimination designs depends sensitively on the alternative under consideration. We finally also note that the $D_2$-optimal design does not require any preliminary information regarding the (unknown) parameters and that the modified $T_{1/6}$-optimal and the $T_{1/3}$-optimal design were constructed for the particular alternative $(c_0, d_0) = (0, 1)$ corresponding to the "true" model. Therefore we expect these designs to be particularly powerful in the examples considered in the simulation study.

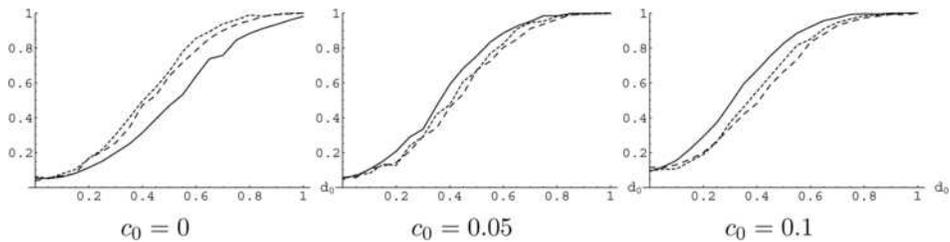

FIG. 2. *Simulated rejection probabilities of the F-test $H_0 : (c_0, d_0) = (0, 0)$ based on the $D_2$-optimal design (dashed line), the modified $T_{1/6}$-optimal design (solid line) and the $T_{1/3}$-optimal design (dotted line) for the parameters $(c_0, d_0) = (0, 1)$ in the cubic regression model (5.2). The errors are centered normally distributed with variance $\sigma^2 = 0.1$.*

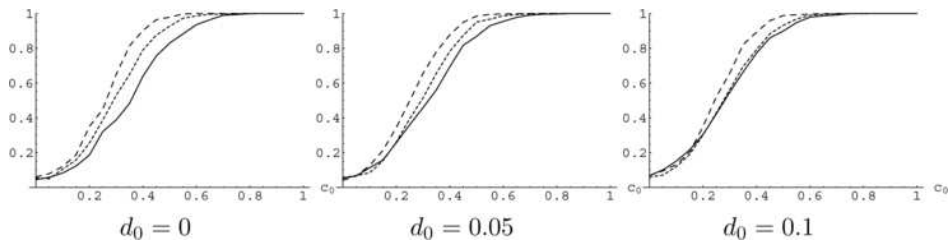

FIG. 3. *Simulated rejection probabilities of the F-test $H_0 : (c_0, d_0) = (0, 0)$ based on the $D_2$-optimal design (dashed line), the modified $T_{1/6}$-optimal design (solid line) and the $T_{1/3}$-optimal design (dotted line) for the parameters $(c_0, d_0) = (0, 1)$ in the cubic regression model (5.2). The errors are centered normally distributed with variance $\sigma^2 = 0.1$*

Usually the next step after model identification is the statistical analysis based on the identified model. Therefore it is also of interest to investigate the performance of the three discrimination designs for this purpose. In Table 3 we present the mean squared errors of the least squares estimates $\hat{a}$, $\hat{b}$, $\hat{c}$ and $\hat{d}$ based on data obtained from a $D_2$-optimal design, the (modified) $T_{1/6}$-optimal and the $T_{1/3}$-optimal design for the special choice of the parameters $(c_0, d_0) = (0, 1)$. The model under consideration is in fact the cubic regression $1 + x + x^3$, for which the $T$-optimal designs were constructed. We observe that the mean squared error of the estimates obtained from the (modified) $T_{1/6}$-optimal design is substantially larger compared to the mean squared error obtained from the $D_2$-optimal and $T_{1/3}$-optimal design. For the last named designs the situation is very similar, where there are slight advantages for the $D_2$-optimal design with respect to the estimation of the parameters $a$ and $c$ and the opposite behavior can be observed for the estimates of the parameters $b$ and $d$.

TABLE 3
*Mean squared error of the least squares estimates in the cubic regression model. The data is obtained from the $D_2$-optimal, the $T_{1/3}$-optimal and (modified) $T_{1/6}$-optimal design for the special choice of the parameters $(c_0, d_0) = (0, 1)$. The variance is chosen as $\sigma^2 = 0.1$*

|  | $D_2$-optimal design | Modified $T_{1/6}$-optimal design | $T_{1/3}$-optimal design |
|---|---|---|---|
| MSE($\hat{a}$) | 0.0050 | 0.0103 | 0.0060 |
| MSE($\hat{b}$) | 0.0290 | 0.0324 | 0.0220 |
| MSE($\hat{c}$) | 0.0120 | 0.0545 | 0.0160 |
| MSE($\hat{d}$) | 0.0360 | 0.0766 | 0.0320 |

5.2. *A nonlinear example.* In this section we consider the problem of discrimination between the exponential regression models

(5.3) $\quad \eta_1(x, \theta_{(1)}) = \theta_{(1)1} \exp(-\theta_{(1)2} x) + \theta_{(1)3} \exp(-\theta_{(1)4} x),$

(5.4) $\quad \eta_2(x, \theta_{(2)}) = \theta_{(2)1} \exp(-\theta_{(2)2} x),$

where the explanatory variable varies in the interval $\mathcal{X} = [-1, 1]$. These models have numerous applications in pharmacokinetics (see, e.g., Shargel and Yu [31] or Rowland [29]) and optimal designs have been discussed extensively in the recent literature (see, Dette, Melas and Pepelysheff [15] or Biedermann, Dette and Pepelysheff [6]). It follows by similar arguments as given in Example 4.3 that a $T$-optimal design has at most three support points. The $T$-optimal designs are listed in Table 4 for various combinations of the parameters $\theta_{(1)j}$, $j = 1, \ldots, 4$.

We have again performed a small simulation study in order to study the rejection probabilities of the likelihood ratio test of the hypothesis

(5.5) $\quad\quad\quad\quad\quad\quad\quad\quad H_0 : \theta_{(1)3} = 0,$

where the data is generated by the different designs. Because this test requires measurements at at least 4 locations, we have modified the $T$-optimal designs by putting 2% of the observations at a fourth point.

For a comparison, there are now two natural candidates based on the $D_s$-optimality criterion. The first design is obtained by maximizing the power of the test for the hypothesis (5.5) and corresponds to the $D_1$-criterion, while the second design is a $D_2$-optimal design in the sense of Stigler [34], and corresponds to the test for the hypothesis

(5.6) $$H_0 : (\theta_{(1)3}, \theta_{(1)4}) = (0, 0).$$

The corresponding local optimal designs are presented in Table 5.

We have simulated data according to the model $\eta_1$ and calculated the power of the likelihood ratio test for the hypothesis (5.5) in various situations. The errors are normally distributed with variance $\sigma^2 = 0.05$ ($\theta_{(1)2} = -1, \theta_{(1)4} = -2$) and $\sigma^2 = 0.2$ ($\theta_{(1)2} = -1, \theta_{(1)4} = 2; \theta_{(1)2} = 2, \theta_{(1)4} = 4$), the sample size is $n = 50$ and 1000 simulation runs are used to calculate the rejection probabilities. Some typical results are depicted in Figure 4, which shows the probability of rejection as a function of the parameter $\theta_{(1)3}$.

TABLE 4
*T-optimal designs for discriminating between the exponential regression models given in (5.3) and (5.4) for a special choice of the parameter $\theta_{(1)}$*

| $\theta_{(1)} = (\theta_{(1)1}, \theta_{(1)2}, \theta_{(1)3}, \theta_{(1)4})$ | $x_1$ | $x_2$ | $x_3$ | $\omega_1$ | $\omega_2$ | $\omega_3$ |
|---|---|---|---|---|---|---|
| $(1, 2, 1, 4)$ | $-1$ | $-0.8$ | $-0.02$ | 0.088 | 0.22 | 0.692 |
| $(1, -1, 1, -2)$ | $-1$ | 0.6 | 1 | 0.645 | 0.246 | 0.109 |
| $(1, -1, 1, 2)$ | $-1$ | $-0.272$ | 1 | 0.168 | 0.437 | 0.395 |
| $(-1, 1, -1, 2)$ | $-1$ | $-0.59$ | 1 | 0.109 | 0.252 | 0.639 |
| $(-1, -1, -1, -0.5)$ | $-1$ | 0.35 | 1 | 0.394 | 0.425 | 0.181 |

TABLE 5
*$D_s$-optimal designs, $s = 1, 2$, for discriminating between the exponential regression models given in (5.3) and (5.4) for a special choice of the parameter $\theta_{(1)}$*

| $(\theta_{(1)1}, \theta_{(1)2}, \theta_{(1)3}, \theta_{(1)4})$ | $s$ | $x_1$ | $x_2$ | $x_3$ | $x_4$ | $\omega_1$ | $\omega_2$ | $\omega_3$ | $\omega_4$ |
|---|---|---|---|---|---|---|---|---|---|
| $(1, 2, 1, 4)$ | 1 | $-1$ | $-0.859$ | $-0.394$ | 0.717 | 0.087 | 0.197 | 0.257 | 0.459 |
| | 2 | $-1$ | $-0.838$ | $-0.404$ | 0.52 | 0.144 | 0.258 | 0.206 | 0.392 |
| $(1, -1, 1, -2)$ | 1 | $-1$ | $-0.03$ | 0.758 | 1 | 0.293 | 0.346 | 0.249 | 0.112 |
| | 2 | $-1$ | 0.03 | 0.697 | 1 | 0.308 | 0.253 | 0.281 | 0.158 |
| $(1, -1, 1, 2)$ | 1 | $-1$ | $-0.636$ | 0.394 | 1 | 0.142 | 0.444 | 0.311 | 0.103 |
| | 2 | $-1$ | $-0.616$ | 0.313 | 1 | 0.341 | 0.309 | 0.268 | 0.082 |
| $(-1, 1, -1, 2)$ | 1 | $-1$ | $-0.758$ | 0.03 | 1 | 0.112 | 0.249 | 0.346 | 0.293 |
| | 2 | $-1$ | $-0.697$ | $-0.03$ | 1 | 0.158 | 0.281 | 0.253 | 0.308 |
| $(-1, -1, -1, -0.5)$ | 1 | $-1$ | $-0.273$ | 0.657 | 1 | 0.215 | 0.631 | 0.29 | 0.134 |
| | 2 | $-1$ | $-0.242$ | 0.576 | 1 | 0.324 | 0.271 | 0.275 | 0.13 |

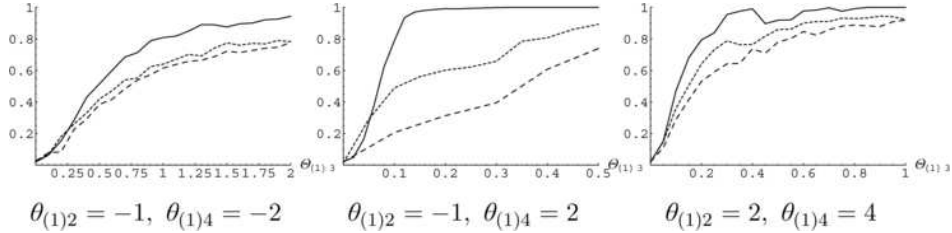

FIG. 4. *Simulated rejection probabilities of the likelihood ratio test for the hypothesis (5.5) based on the $D_1$-optimal design (dashed line), $D_2$-optimal design (dotted line) and the $T$-optimal design (solid line) in the exponential regression model (5.3), where $\theta_{(1)1} = \theta_{(1)3} = 1$.*

If both parameters in the exponential functions are negative (left panel in Figure 4) the power of the test obtained from the modified $T$-optimal design is larger than the power of the test based on the $D_2$-optimal design. On the other hand the $D_2$-optimal design seems to have slightly better discrimination properties than the $D_1$-optimal design in this case. If both parameters are of opposite sign (middle panel in Figure 4) the situation is different and the $D_2$-optimal design yields a bit more power for small values of the parameter $\theta_{(1)3}$. In this example the $D_1$-optimal design is totally defective. Finally, the right panel of Figure 4 shows a situation where both parameters in the exponential functions are positive. If both parameters are of opposite sign (middle panel in Figure 4) the situation is different and the $D_2$-optimal design yields a bit more power for small values of the parameter $\theta_{(1)3}$. In this example the $D_1$-optimal design is totally defective. Finally, the right panel of Figure 4 shows a situation where both parameters in the exponential functions are positive. Here almost the same behavior as in the case of negative parameters is observed. While the $D_2$-optimal design yields more power than the $D_1$-optimal design, the test based on the (modified) $T$-optimal shows the best performance. On the other hand the $D_2$-optimal design advices the experimenter to take observations at 4 different locations and therefore it also allows the estimation of all parameters in the extended model.

The impact of the discriminating designs on the parameter estimates is investigated in Table 6, where we exemplarily show two typical examples of the simulated mean squared error of the parameter estimates under the different designs. If $\theta_{(A)} = (1, -1, 1, 2)$ the $D_1$- and $D_2$-optimal designs yield substantially smaller mean squared errors than the $T$-optimal design, and the $D_1$-optimal design shows a slightly better performance than the $D_2$-optimal design. In the case $\theta_{(B)} = (1, 2, 1, 4)$ the $D_1$- and $D_2$-optimal design yield the smallest mean squared errors, while the (modified) $T$-optimal shows again the worst performance. The mean squared errors obtained by the $D_2$-optimal design are slightly larger than those obtained by the $D_1$-optimal

Table 6

*Simulated mean squared error of the least squares estimates in the exponential regression model. The data is obtained from the $D_1$-, $D_2$- and (modified) $T$-optimal design for the special choice of the parameters $\theta_{(A)} = (1, -1, 1, 2)$ and $\theta_{(B)} = (1, 2, 1, 4)$. The variance is chosen as $\sigma^2 = 0.2$*

|  |  | $D_1$-optimal design | $D_2$-optimal design | $T$-optimal design |
|---|---|---|---|---|
| $\theta_{(A)}$ | MSE($\hat{a}$) | 0.04491 | 0.05507 | 0.31266 |
|  | MSE($\hat{b}$) | 0.05468 | 0.06687 | 1.07216 |
|  | MSE($\hat{c}$) | 0.02503 | 0.03137 | 0.15910 |
|  | MSE($\hat{d}$) | 0.02414 | 0.02803 | 0.15603 |
| $\theta_{(B)}$ | MSE($\hat{a}$) | 0.18217 | 0.18552 | 0.37235 |
|  | MSE($\hat{b}$) | 0.57880 | 0.80178 | 2.94709 |
|  | MSE($\hat{c}$) | 0.18374 | 0.17361 | 0.37019 |
|  | MSE($\hat{d}$) | 0.25151 | 0.21136 | 0.43687 |

design. Summarizing these and similar results (which are not shown for the sake of brevity) we conclude that the $D_1$- and $D_2$-optimal designs have good properties for model discrimination and additionally have good properties for parameter estimation if the null hypothesis (5.6) has been rejected. In many cases the mean squared error of the parameter estimates obtained from the modified $T$-optimal design is at least two times larger compared to the results obtained from the $D_1$- and $D_2$-optimal designs.

**Acknowledgments.** The authors are grateful to Martina Stein who typed parts of this paper with considerable technical expertise and to Christine Kiss for computational assistance. The authors would also like to thank two unknown referees for their constructive comments on an earlier version of this paper.

RUHR-UNIVERSITÄT BOCHUM
FAKULTÄT FÜR MATHEMATIK
44780 BOCHUM, GERMANY
E-MAIL: holger.dette@rub.de
stefanie.titoff@rub.de